\documentclass[12pt,twoside,leqno]{amsart}  

\usepackage{amsmath,amscd,amssymb,amsfonts,epsfig}

\newcommand{\sO}{{\mathcal O}}

\newcommand{\F}{{\mathbb F}}
\newcommand{\G}{{\mathbb G}}

\newcommand{\K}{{\mathbb K}}

\newcommand{\N}{{\mathbb N}}
\newcommand{\pP}{{\mathbb P}}

\newcommand{\bo}{{\bf B} }
\newcommand{\g}{{\bf G} }

\newcommand{\lra}{\longrightarrow }

\newcommand{\ddim}{{\rm dim}}

\newcommand{\Ker}{{\rm Ker}}

\newcommand{\im}{{\rm im}}

\theoremstyle{plain}
\newtheorem{thm}{Theorem}
\newtheorem{lem}[thm]{Lemma}
\newtheorem{cor}[thm]{Corollary}

\newtheorem{remark}[thm]{Remark}

\newtheorem{defn}[thm]{Definition}

\numberwithin{thm}{section}
\numberwithin{equation}{section}
\def\pf{\noindent{\sc Proof. }}
\def\qed{\rm Q.E.D.}

\begin{document}

\title{On the motive of certain subvarieties of fixed flags}

\author{Pedro Luis del Angel}
\address{Pedro Luis del Angel}
\address{\vskip-6.5mm Cimat, Guanajuato, Mexico}
\email{luis@cimat.mx}

\author{E. Javier Elizondo }
\address{E. Javier Elizondo }
\address{\vskip-6.5mm INSTITUTO DE MATEMATICAS \\
  Universidad Nacional Aut\'onoma de M\'exico \\ Ciudad Universitaria \\
 04510, Mexico DF.}
\email{javier@math.unam.mx}

\date{}

\begin{abstract} We compute the Chow motive of certain subvarieties of the
Flags manifold and show that it is an Artin motive.
\end{abstract}
\thanks{Partially supported  by DGAPA IN105905}

\maketitle

\section*{Introduction}

Let $\g$ be a connected algebraic semisimple group defined over an
algebraically closed field $\K$, with universal separable covering. We
denote by $U$ the variety of unipotent elements of $\g$ and by $B$ a
Borel subgroup of $\g$. It is well known  (see \cite{steinberg}) that
if $\bo = \g/B$ and 
$$
Y:=\{(x,gB)\in U\times \bo | g^{-1}xg \in B\}
$$
then
$$
\pi: Y\lra U
$$
is a desingularization,  $\pi$ denotes the natural projection.\\

 Consider the case   $\g=SL_n$, then it is easy to see that  
the variety  of complete flags is $\bo=\F$. For this case we have that for any
unipotent element $x$ 
the fiber $\pi^{-1}(x)$ is isomorphic to the 
variety of fixed flags. Moreover, J.A. Vargas
(see \cite{vargas}) has given a description of a dense open set for
every irreducible component of the fiber, and N. Spaltenstein (see
\cite{spaltenstein}) has constructed a stratification of the fiber,
which unfortunately is not completely compatible with the
decomposition into the irreducible components.\\ 

The purpose of this work is to describe the motive of the
irreducible components of the fibers $\pi^{-1}(x)$ when $x$ is of
type $(p,q)$. Observe that it is interesting to know the geometry 
and K-theory of the fibers of the  desingularization since algebraic
singularities appear  
within the unipotent variety. It is also important
for applications such as computing zeta functions and counting points
over finite fields.\\
 
 The paper is divided as follows. In the first and second section we
introduce  basic
notation and the stratification given by  Spaltenstein. In the third 
section it is  given a description of the irreducible components of the
fiber of $x$ for  any unipotent element $x$ of type $(p,q)$. In the fourth
section  we use the above description to compute the Chow motive,
and show that the image of the Chow motive of any irreducible
component into the Voevodski's category is an Artin motive. We also compute the
motive of some irreducible components of a slightly more general type, showing
that these motives are extension by Artin motives of the motive of a product
of flag varieties.\\

 The authors want to thank Pedro Dos Santos, Bruno Kahn, Jochen
Heinloth, James Lewis  
and Stefan M\"uller-Stach for some useful discussions and
suggestions. We are particularly  
indebted to Herbert Kurke for pointing out a mistake in a previous
version of this article. 

\section{Preliminaries}

Let $\g$ be the group $SL_n$ with coefficients in an algebraically
closed field $\K$. 
Consider  a Borel subgroup $B$ of $\g$ and $T$ a maximal torus on
$B$. If $V$ is a $\K$-vector space of dimension $n$, then the
variety $\bo$ is isomorphic to the variety of complete flags
$\F=\F(V)$. \\
For any unipotent element $x\in U\subset \g$, we denote the fiber 
$\pi^{-1}(x)$ by $\F_x$ . One says that $x$ is of type 
$(\lambda_1, \cdots, \lambda_s)$ if the 
Jordan canonical form of $x$ consists exactly of $s$ blocks of sizes
$\lambda_1\ge\lambda_2\ge\cdots\ge\lambda_s$.

If we write $x=1+n$ where $n$ is the nilpotent part of $x$, then there
is a basis 
$$
\{ e_{i,j}| 1\le j\le s\; , 1\le i\le \lambda_j \}
$$
of $V$ adapted to $x$ in the sense that\,  $n(e_{i,j})=e_{i-1,j}$, where
$e_{0,j}=0$\, for every $j$. Therefore we can write
$V=V_1\oplus \cdots\oplus V_{\lambda_1}$ and
$n:V_{\lambda_1}\lra V_{\lambda_1-1}\lra \cdots V_1\lra 0$, where
$V_i$ is the space generated by $\{ e_{i,j}\}$ with $i$ fix.

It is well known that if $x\in U\subset SL_n$ is unipotent of type
$(\lambda_1, \cdots, \lambda_s)$, then the fiber $\F_x$ has as many
irreducible components as standard tableaux of type
$(\lambda_1, \cdots, \lambda_s)$. The shape of a standard tableau is as follows:
\begin{figure}[h!]
        \begin{center}
         \scalebox{.20}{\includegraphics{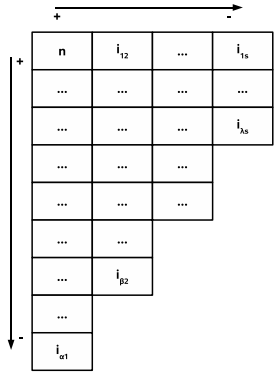}}
        \end{center}
        \caption{Standard tableau of type $(\alpha,\beta,\cdots,\lambda)$.}
\end{figure}

\noindent where the numeration strictly decrease from the top to the bottom
and from left to right.
Spaltenstein constructed a stratification of $\F_x$. We will use,
and explain, the particular case where $x$ is of type $(p,q)$.

Given $x\in U\subset SL_n$ of type $(p,q)$ and a $\K$-vector space $V$ of
dimension $n$, consider a basis
$\{ e_{i,j}| 1\le j\le 2\; , 1\le i\le \lambda_j \}$
of $V$ adapted to $x$, as mentioned above. Given a number $t\in\{1,2\}$ 
we define some 
subsets of $\F_x$ as follows:
\begin{enumerate}
\item For $t=1$, let $X_1=X_1(V)$ be the set of flags
$F_1\subset\cdots\subset F_n \in \F$ such that \\$F_1:= <e_{1,1}>$.
\item For $t=2$, let $X_2=X_2(V)$ be the set of flags
$F_1\subset\cdots\subset F_n \in \F$ such that \\$F_1:= <a e_{1,1}+
e_{1,2}>$ for some number $a\in \K$ not necessarily different from
zero.
\end{enumerate}
We also  define inductively the sets
$$
X_{i j}(V):=\{F_1\subset\cdots\subset F_n\in X_j(V)\; |\;
F_2/F_1\subset\cdots\subset F_n/F_1 \in X_i(V/F_1)\}.
$$
Applying the same process to $X_{i j}(V)$, we get sets of the form
$X_{i j k}(V)$ and
so on. It is
not difficult to see that after $n$ times we end up with a locally
closed subset of flags which actually belongs to 
$\F_x$, since a flag $F:= F_1\subset\cdots\subset F_n$ is in $\F_x$
if and only if $n(F_i)\subset F_{i-1}$ for every $i$. Moreover the
sets we obtain  form a stratification of $\F_x$.

It is clear from the construction that all spaces of
the stratification are affine spaces and that the
affine strata of maximal dimension
(which is precisely $q$ if $x$ is of type $(p,q)$) are open disjoint sets. 
You can also count the number of irreducible
components of $\F_x$, which coincides with the number of standard tableaux 
of the given type.

\begin{remark}
The natural projection from $X_1(V)$ to $\F(V/F_1)$ is an isomorphism.
\label{isomorphism}
\end{remark}

\begin{remark}
\label{component}
Given a nilpotent element $x\in U$ of type 
$(\lambda_1, \cdots, \lambda_s)$,
a standard tableau $\sigma$ as above, a
$\K$-vector space $V$ of dimension $n=\lambda_1+\cdots\lambda_s$ and a
basis $\{e_{i_1,\cdots,i_s}\}$ of $V$ adapted to $x$, there is a
maximal affine space 
among those obtained from Spaltenstein's stratification of $\F_x$,
which can be associated to $\sigma$.

Let $\psi:\{1,\cdots,n\}\lra\{1,\cdots,s\}$ be the function defined as
$\psi(k)= i$, if $k$ appears in the $i-$th column of $\sigma$,
counting from left to right. Now, let us 
consider the stratum $X_{\sigma}:=X_{\psi(n)\cdots\psi(1)}$, then $X_{\sigma}$
is of maximal dimension and 
the corresponding irreducible component will be denoted by $Y_{\sigma}$.
\end{remark}

\section{Decomposable irreducible components}

Let $A:=\{a_1,\cdots, a_r\}$ be a totally ordered set and we assume
that $a_1<\cdots <a_r$.  
Let $\phi:A \lra \{1,\cdots,s\}$ be a surjective map. Then the pair
$(A,\phi)$ induces a tableau 
$\sigma:=\sigma(A,\phi)$ as follows.

The tableau $\sigma$ has $s$ columns, the $i$-th column of the tableau has 
$b_i:=|\phi^{-1}(i)|$ boxes and one fills them up according to $\phi$,
i.e. starting with 
$r$ and then going down to $1$. 
The number $k$ should appear in the column $\phi(k)$. The numbers
appearing in a 
column are in decreasing order, and the tableau will be of
type $(b_1,\cdots,b_s)$. 

Not every tableau obtained in this way is a Young tableau,
unless $b_1\ge\cdots\ge b_s$, but it might be no standard. For
instance, the tableau associated to  
$\phi(1)=\phi(2)=1$, $\phi(3)=\phi(4)=2$
and $\phi(5)=3$ is a Young tableau but it is not standard,
whereas the tableau 
associated to $\phi(1)=\phi(4)=2$, $\phi(2)=3$ and $\phi(5)=\phi(3)=1$
is a standard Young tableau. 

If $\phi$ is a decreasing bijection and $|A|>1$ we say that the pair $(A,\phi)$ 
is of flag type.
\begin{remark}
Any standard Young tableau of type $(\lambda_1\cdots \lambda_s)$ is
obtained  
from the pair $(\{1,\cdots ,n\}, \psi)$, where $\psi:\{1,\dots,n\}\lra
\{1,\cdots,s\}$ is the  
function described in \ref{component}.
\end{remark}

\begin{defn}
A standard Young tableau $\sigma$ of type $(\lambda_1\cdots \lambda_s)$, with 
$n=\lambda_1 + \cdots + \lambda_s$, is called {\em decomposable} if
there is a partition 
${\displaystyle \{1,\cdots ,n\}=(\sqcup_k A_k)\bigsqcup (\sqcup_t B_t)}$
by totally ordered sets 
such that the pair
$(A_k,\psi|_{A_k})$ is of flag type for all $k$, i.e. if $|A_k|=m_k$ then
$\psi|_{A_k}$ is a decreasing bijection between $A_k$ and $\{1,\cdots,m_k\}$, 
and the pair $(B_t,\psi|_{B_t})$ induces
a standard Young tableau of type $(p_t,q_t)$ for every $t$, where $\psi$
is as above. In particular \,
$\im \,(\psi|_{B_t}) \;= \{1,2\}$.
\end{defn}

\begin{remark}
A pair $(A,\phi)$ is of flag type if and only if the associated standard tableau
$\sigma$ is of type $(1,\cdots,1)$,
which corresponds to the identity in $SL_{n}$. In this case there is only one
irreducible component, namely the whole flag variety. 
\end{remark}

\section{Irreducible components of $\F_x$ for
$x$ of type $(p,q)$}

The description of the irreducible components of type
$(p,q)$ given here is already   
in \cite{thesis}. However, the present proof is better than the
one that appears there, and it is also important to incorporate the
description here for the sake of completeness. 

Let $x=1+n$ be a unipotent element of type $(p,q)$ in $SL_n$ with
nilpotent part $n$, 
$V$ a $\K$-vector space of dimension $n=p+q$ and $\F=\F(V)$ the
variety of complete flags on $V$. 

We have the following related three lemmas. 
\begin{lem}
The set \label{cerrado1}
$$
A=\{ F:= F_1\subset\cdots\subset F_n \in\F \; |\; n^i(F_k)\subset F_m\}
$$
is closed in $\F$ for all fixed $i,k,m\in\N$. We set $F_0=0$.
\end{lem}
\begin{lem}
The set \label{cerrado2}
$$
L=\{ F:= F_1\subset\cdots\subset F_n \in\F \; |\; n^i(F_k)\subset S\}
$$
is closed in $\F$ for all fixed $i,k\in\N$. $S$ is a fixed
subspace of $V$. 
\end{lem}
\begin{lem}
The set \label{cerrado3}
$$
H=\{ F:= F_1\subset\cdots\subset F_n \in\F \; |\; dim (F_r+ n^k(F_r))\le d \}
$$
is closed in $\F$ for all fixed $r,k,d\in\N$.
\end{lem}

For these lemmas one shows that the corresponding set
(either $A$, $L$ or $H$) is a determinantal set and therefore algebraic.

Let $\pi_r:\F \lra \G r(1,n)\times\cdots\times\G r(r,n)$ be the
composition of the natural embedding of $\F$ in 
$\G r(1,n)\times\cdots\times\G r(n-1,n)$ followed by the projection to the
first $r$ factors.
\begin{thm}
Let $\sigma$ be a standard tableau of type $(p,q)$, and 
$Y_{\sigma}$  the corresponding irreducible component of $\F_x$. Let
$1\le r\le n$ be 
a natural number and set 
$Y_{\sigma}(r):=\pi_r(Y_{\sigma})$. We set $Y_{\sigma}(0)$ to be a
fixed point. 
Let $f_r$ be the natural projection
$Y_{\sigma}(r)\lra Y_{\sigma}(r-1)$,
then for all $p\in Y_{\sigma}(r-1)$ one has
$$
f_r^{-1}(p)=\left\{\begin{matrix}
                         1pt  & \mbox{if} & r & \mbox{appears in the
                           left column of} & \sigma,  \cr
                         \pP^1 & \mbox{if} & r & \mbox{appears in the
                           right column of} & \sigma. 
                  \end{matrix}\right.
$$
\ \\
In particular, $Y_{\sigma}(r)\lra Y_{\sigma}(r-1)$ is either a
$\pP^1$-bundle or an isomorphism for 
every $1\le r\le n$. \label{irreducible}
\end{thm}

\pf 
Let us assume that there exist
a 2-dimensional hermitian $\K$-Vector space $W$ (see II below). The more general
situation follows in a similar way, but it does not give more light into
the geometry of the problem, which is our interest.

As mentioned before, there is an open subset $X_{\sigma}$ of the
irreducible component $Y_{\sigma}$.
Consider the image $X_{\sigma}(r):=\pi_r(X_{\sigma})$ of this set under the map
$\pi_r$. It is not hard to observe that the fibers of $f_r$ restricted to
$X_{\sigma}(r)$ are isomorphic to a point 
or to an affine line, depending on whether $r$ appears in the left or in the 
right column of $\sigma$. Now, we need to consider the map $f_r$ as a map
from $Y_{\sigma}(r)$, and for this we need to give a good
description of $Y_{\sigma}$. We will do these in four
steps, which involve different cases. The three last lemmas will be used. 
\ \newline

\noindent {\bf I.} \hspace{.3cm}Fix a basis $\{e_{i,j}\}$ of $V$
adapted to $x$ as before, and 
let $\sigma$ be a standard Young tableau of type $(p,q)$. The map 
$\psi:\{1,\cdots ,n\}\to \{1,2\}$ defined in remark \ref{component}. 
induces a partition 
${\displaystyle \{1,\cdots ,n\}=(\sqcup_{k=1}^t A_k)\; \bigsqcup\;
  (\sqcup_{k=1}^t B_k)}$ 
such that
$A_k\subset\psi^{-1}(1)$ and $B_k\subset\psi^{-1}(2)$ are made up by consecutive 
integers for all $k$. For $1\le i\le t$ define $a_i=|A_i|$, $b_i=|B_i|$, 
$s_j=\sum_{i<j}\; a_i$ and $S_j =\sum_{i<j}\; b_i$.
We also set $a_0=0$ and $b_0=0$.

Because $\sigma$ is a standard tableau, it is
clear that either  
$A_1<B_1<A_2<B_2<\cdots <A_t$ and $B_t=\emptyset$ or $B_1<A_1<B_2<A_2<\cdots <B_t<A_t$,
where $C<D$ means $c<d$ for all $c\in C$ and all $d\in D$.

Without loss of generality we can assume $B_1<A_1<\cdots <B_t<A_t$. Indeed,
if $A_1<B_1$ then $A_1=\{1, \cdots , a_1\}$, in this case $p-a_1>q$
since $\sigma $ is a standard tableau and the numeration inside of it
decreases to the right.
Therefore, for every flag\, 
$F_1\subset \cdots \subset F_n\in X_{\sigma}$ \, one has\ $F_i=\im \;
n^{p-i}$\, for\, $1\le i\le a_1$. 
But these are closed conditions and so they are fulfilled by
all flags in $Y_{\sigma}$. In this situation we have that remark
\ref{isomorphism} implies that\,  
$Y_{\sigma}\cong Y_{\sigma'}(V/(\im\;n^{p-a_1}))$; where $\sigma'$ is the
standard Young tableau induced by \,
$\left(\{a_1+1,\cdots ,n\},\psi|_{\{a_1+1,\cdots ,n\}}\right)$.  
\\ \ \\
{\bf II.}\hspace{.3cm} 
We first consider the following case.
\begin{figure}[h!]
        \begin{center}
         \scalebox{.15}{\includegraphics{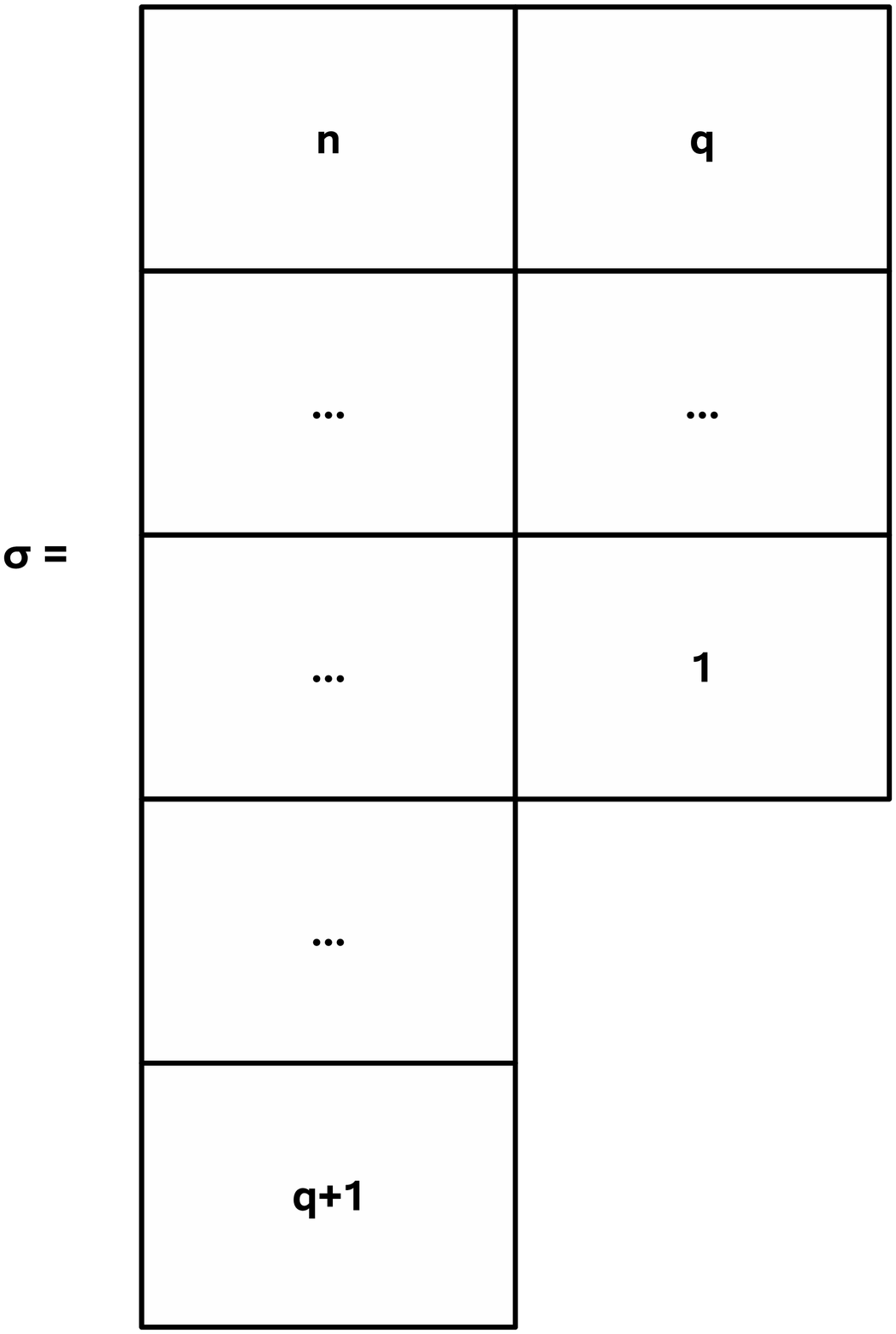}}
        \end{center}
\end{figure}
\newline
 Let $W$ be an hermitian vector space of dimension $2$ with basis
 $\{w_1,w_2\}$. For every point 
$(P_1,\ldots, P_q)\in \pP(W)^q$ write $P_i=(a_i:b_i)$.
 Let $R_i=(c_i:d_i)$ be the point in $\pP(W)$
that represents the orthogonal complement of the subspace\,
$<a_iw_1+b_iw_2>\,\subset W$.  
\\ 
For all $1\le j\le q$ consider the vectors (in $V$)
$$
\begin{matrix}
P_{1,j} & = & a_1e_{j,1}+b_1e_{j,2} \cr
R_{1,j} & = & c_1e_{j,1}+d_1e_{j,2}
\end{matrix}
$$
and for all $i,j$ with $i>1$ and $i+j\le q+1$ define the following vectors
$$
\begin{matrix}
P_{i,j} & = & a_iR_{i-1,j}+b_iP_{i-1,j+1} \cr
R_{i,j} & = & c_iR_{i-1,j}+d_iP_{i-1,j+1}
\end{matrix}
$$
Now, with the notation described above, and for $1\le k\le q$, let $Y_k$
 be the set  giving by the following conditions.
\\ $F_1\subset\cdots\subset F_n\in Y_k$ \, if and only if
there is a point  \, $(P_1,\ldots ,P_q)\in\pP(W)^k\times\left[ \pP(W)\; -\;
  \{(1:0)\}\right]^{q-k}$ 
such that
\begin{enumerate}
\item  $F_m = <P_{1,1},\ldots,P_{m,1}>$ \, for all\, $1\le m\le q$.
\item  $F_{q+t} = <P_{1,1},\ldots,P_{q,1},R_{q,1},\ldots,R_{q-t+1,t}>$.
        \,for all\, $1\le t\le q$ 
\item $F_{s} = \Ker\; n^s$ \, for all \, $2q\le t\le n$.
\end{enumerate}
Observe that if the $P_i$'s are different from $(1:0)$ for all $i$, then
$$<P_{1,1},\ldots,P_{q,1},R_{q,1},\ldots,R_{q-t+1,t}>=
<P_{1,1},\ldots,P_{q,1},e_{1,1},\ldots,e_{1,t}>$$
 for all $t$, therefore
$Y_0=X_{\sigma}$. 
This situation corresponds to:\newline
$F_1\subset\cdots\subset F_n\in X_{\sigma}$ if and only if
\begin{enumerate}
\item  $F_1\subset \Ker\; n \; - \; <e_{1,1}>$,
\item $n(F_{i})\subset F_{i-1}$ but $n(F_{i})\not\subset F_{i-2}$ for every
      $2\le i\le q$,
\item $\ker \; n^t \subset F_{q+t}$ for every $ 1\le t\le p$.
\end{enumerate}
Similarly the sets $Y_k$ for $1\le k<q$ satisfy
$F_1\subset\cdots\subset F_n\in Y_k$  
if and only if these conditions are fulfilled:
\begin{enumerate}
\item $F_m\subset \Ker\; n^m$ for all $1\le m\le k$
\item $n(F_{i})\subset F_{i-1}$ but $n(F_{i})\not\subset F_{i-2}$ for every
       $k+1\le i\le q$
\item $\ker\; n^t\subset F_{q+t}$ for every $ 1\le t\le p$.
\end{enumerate}
Finally the set $Y_q$ satisfies $F_1\subset\cdots\subset F_n\in Y_k$ 
\, if and only if the three following conditions are fulfilled:
\begin{enumerate}
\item $F_m\subset \Ker\; n^m$ for all $1\le m\le k$,
\item $n(F_{i})\subset F_{i-1}$, for all $1\le i\le n$,
\item $\ker\; n^t\subset F_{q+t}$ for every $ 1\le t\le p$.
\end{enumerate}
where we have set $F_0=0$.

Now $Y_q$ contains $X_{\sigma}=Y_0$ and is irreducible by construction, 
moreover it is actually a closed set because of lemmas \ref{cerrado1},
\ref{cerrado2} and 
\ref{cerrado3}. Therefore $Y_q=Y_{\sigma}$. 
\\
The natural projection $f_{k}:Y_{\sigma}(k)\lra Y_{\sigma}(k-1)$ 
is nothing more than the map
$$
F_0\subset\cdots\subset F_{k-1}\subset F_{k}\mapsto
F_0\subset\cdots\subset F_{k-1} 
$$
and, because of the construction of $Y_q$, the fiber of this map is
$\pP^1$ for all $1\le k\le q$. On the other hand
\begin{eqnarray}
\label{fiber}
<P_{1,1},\cdots,P_{q,1},R_{q,1},\cdots,R_{q-k+1,k}>=
<P_{1,k+1},\cdots,P_{q-k,k+1}>\oplus\mbox{ ker }n^k
\end{eqnarray}
for all $1\le k\le q$, therefore $f_{q+k}$ is an isomorphism for $1\le
k\le q$ since the 
space $F_{q+k}$ is already determined by the space
$F_{q+k-1}$. Finally, we have that  
$F_{2q+t}=\ker\; n^{q+t}$  for all $t\ge 0$, and the theorem follows
in this case. 
\\ \\
{\bf III.}\hspace{.3cm} Let us assume that $s_j< S_j$ for all
$1\le j < t$.
We follow the notation introduced in  [I] and [II]. We define $Y_k$ in a
similar way than in the other cases. Namely, \newline
$F_1\subset\cdots\subset F_n=V\in Y_k$\, if and only if there exist a
point \,  
$(P_1,\cdots ,P_q)\in\pP(W)^k\times\left( \pP(W)\; -\; \{(1:0)\}\right)^{q-k}$
such that
\begin{enumerate}
\item $F_{S_j + m} = \ker\; n^{s_j}\;\oplus <P_{1,s_j+1}, \cdots
  ,P_{S_j-s_j+m,s_j+1}>$ \, for all\, $1\le m\le b_j$, and \, $1\le j\le t$.
\item $F_{S_{j+1}+s_j+m} = \ker\; n^{s_j}\;\oplus
  <P_{1,s_j+1},\cdots,P_{S_{j+1}-s_j,s_j+1},R_{S_{j+1}-s_j,s_j+1},
  \cdots,R_{S_{j+1}-s_j-m+1,s_j+m}>$ \, for all\, $1\le m\le a_j$, and
  $1\le j\le t$. 
\item $F_{s} = \Ker\; n^s$ for all $2q\le t\le n$
\end{enumerate}
One sees immediately that $X_{\sigma}=Y_0\subset \cdots \subset
Y_q$. Observe that 
$Y_q$ is irreducible  by construction, and similarly as in [I], it can
be described by 
the following conditions\newline
$F_1\subset\cdots\subset F_n\in Y_k$ \
if and only if
\begin{enumerate}
\item $\ker\; n^{s_j}\subset F_{S_j+m}\subset \Ker\; n^{S_j+m}$ for all
     $1\le m\le b_j$, for all $1\le j\le t$,
\item $\ker\; n^{s_j+m}\subset F_{S_{j+1}+s_j+m}\subset \ker\; n^{S_{j+1}+m}$ 
for all $1\le m\le a_j$, for all $1\le j\le t$,
\item $\ker\; n^{m}=F_{m}$ for all $2q\le m\le p$.
\end{enumerate}
where condition (2) is a consequence of the following equality
\begin{eqnarray}
\begin{matrix} \label{fiber2}
\mbox{ker }n^{s} &  \oplus
<P_{1,s+1},\cdots,P_{m,s+1},R_{m,s+1},\cdots,R_{m-t+1,s+t}>\cr  & \cr
 & = \mbox{ ker }n^{s+t}\oplus <P_{1,s+t+1},\cdots,P_{m-t,s+t+1}>.
\end{matrix}
\end{eqnarray}
where $1\le t\le m$ and $s+t+1\le q$.
\\
As in [II], it follows from the construction of $Y_q$ that the fiber
of $f_k$ is isomorphic to $\pP^1$ if $k\in B_s$, for some $s$. We have
that either from
condition (3) or from equation \ref{fiber2} above, one gets that
$f_k$ will be an isomorphism if $k\in A_s$ for some s.
\\ \\
{\bf IV.}\hspace{.3cm} We follow the notation introduced in
[III]. Consider a standard tableau $\sigma$ such that 
$s_j\not< S_j$ for some $1\le j\le t$. Let $j_0$ be the smallest index 
such that
$s_j\ge S_j$. Since the $A_i$'s and the $B_i$'s are made up by
consecutive numbers, 
and for every index $k>m$ the numbers appearing in \,$B_k$\, and
\,$A_k$\, are bigger 
than those appearing in \,$A_m$\, or\, $B_m$\, by applying induction one shows
that $\max\{A_m\}=s_{m+1}+S_{m+1}$\, and\, $\max\{B_m\}=s_m+S_{m+1}$
for every $m$, in particular $\max\{A_{j_0}\}\ge 2S_{j_0}$. Moreover,
since $j_0$ was 
minimal then $\max\{B_{j_0}\}< 2S_{j_0}$, i.e. $2S_{j_0}\in
A_{j_0}$. Therefore, for 
every flag\, $0=F_0\subset\cdots\subset F_n\in X_{\sigma}$\, one has 
$F_{2S_{j_0}}=\ker n^{S_{j_0}}$. Since this is a closed condition, this
is also true  
for all flags in the irreducible component $Y_{\sigma}$, and therefore
$Y_{\sigma}\cong Y_{\sigma'}(\ker\; n^{S_{j_0}})\times
Y_{\sigma''}(V/\ker\; n^{S_{j_0}})$,  
where $\sigma'$ is the standard 
tableau induced by the pair $([1,\cdots ,2S_{j_0}]; \psi|_{[1,\cdots ,2S_{j_0}]})$,
and $\sigma ''$ is the standard tableau induced by the pair
$([2S_{j_0}+1,\cdots , n]; \psi|_{[2S_{j_0}+1,\cdots ,n]})$, both of them of type
$(a,b)$ for some $a\ge b$.
The proposition follows by induction on the dimension.

\hskip1cm{}\hfill \qed

\section{The motive of the irreducible
components}

If $X$ is a scheme and $G$ a group, the $G$-torsors on $X$ for the \'etale
cohomology are parametrized by
$H^1_{et}(X,G)$. The exact sequence
\begin{eqnarray*}
0\lra \G_m\lra GL_2\lra PGL_2\lra 0
\end{eqnarray*}
gives  a map
$$
\delta:\; H^1_{et}(X, PGL_2)\lra H^2_{et}(X,\G_m)=Br(X).
$$
Moreover, $z\in\mbox{ ker }\delta \Leftrightarrow z$ can be extended
to a $GL_2$ torsor on $X$, i.e. if $z$ can be extended to a
vector bundle on $X$ for the \'etale topology. Therefore, in order
to show that a torsor that corresponds to an element $z\in H^1_{et}(X,
PGL_2)$ is the projective bundle associated to a rank $2$ vector
bundle on $X$, it is enough to show that its image in $Br(X)$ is zero.
In this chapter we will deal with varieties over an algebraically 
closed field $\K$ of characteristic zero and therefore the Brauer 
group of $X$ coincide with the ``geometric'' Brauer group for this case.

\begin{thm}
Let $\sigma$ be a standard tableau
of type $(p,q)$ and $Y_{\sigma}$ be the corresponding irreducible component of
$\F_x$. Then the motive $h(Y_{\sigma})$ is isomorphic to $(1+L)^q$. In
particular, it is an Artin motive.
\end{thm}

\pf Since the open cell $X_{\psi(1), \cdots,\psi(n)}\subset Y_{\sigma}$ is an 
affine space, it follows from the proof of theorem \ref{irreducible}
that all the 
varieties $Y_{\sigma}(r)$ are rational and so
$Br(Y_{\sigma}(r))=0$ for all $r$, therefore
$Y_{\sigma}(r+1)\lra Y_{\sigma}(r)$ is either an isomorphism or
the projective bundle associated to a
rank $2$ vector bundle on $Y_{\sigma}(r)$, in this case we have 
$h(Y_{\sigma}(r+1))\cong (1+L)\otimes h(Y_{\sigma}(r))$, see \cite{manin}.
Since $\ddim \; Y_{\sigma}=q$ the conclusion follows by induction.

\hskip1cm{}\hfill\qed

\begin{remark}
Following the notation of remark \ref{component}, if $Y_{\sigma}$ is an
irreducible component for which 
$$\psi (t)=\left\{ 
\begin{matrix}2 & \mbox{for $t$ odd} \cr
 & \cr
1 & \mbox{for $t$ even}
\end{matrix}
\right.
$$
then $ Y_{\sigma}\cong (\pP^1)^q $, and the multiplicative structure of the 
corresponding motive is clear. In general one needs to find sections of the 
maps $Y_{\sigma}(r+1)\to Y_{\sigma}(r)$ that correspond to the bundle
$\sO_{Y_{\sigma}(r+1)}(1)$ and compute their autointersection 
numbers to explicitly find a normalized rank 2 vector bundle which induces the 
$\pP^1$- bundle over $Y_{\sigma}(r)$ and, lastly, being able to actually 
compute the multiplicative structure of the motive, see \cite{manin}.
\end{remark}
\begin{cor}
If $\sigma$ is a decomposable standard tableau then its motive is an Artin motive.
\end{cor}

\pf If $\sigma$ is decomposable then the irreducible component $Y_{\sigma}$
is isomorphic to a product of towers of $\pP^1$-bundles over flag
varieties. Since flag  
varieties are rational, then all the Brauer groups involved are
zero. Moreover, since 
the flag varieties are themselves towers of projective bundles
associated to vector bundles, the 
motive that is obtained is an Artin motive.
\hskip1cm{}\hfill \qed

\bibliographystyle{plain}
\renewcommand\refname{References}

\end{document}